\magnification = 1120
\showboxdepth=0 \showboxbreadth=0
 %no black boxes

\hsize = 5.5truein
\vsize = 8.5truein
\hfuzz = 5pt
\baselineskip 18pt
\parskip3pt
%\voffset = -0.5in
%\hoffset = -.05truein
\def\qed{\hfill\vrule height6pt  width6pt  depth0pt }

\def\ss{\smallskip}
\def\ms{\medskip}
\def\bs{\bigskip}

\def\cl{\centerline}

\def\nind{\noindent}

\def\ref#1#2{\nind\hangindent.5in\hbox to .5in{#1\hfill}#2}
\def\reff#1#2{\nind\hangindent.8in\hbox to .8in{\bf #1\hfill}#2\par}
\def\refd#1#2{\nind\hangindent.8in\hbox to .8in{\bf #1\hfill {\rm
--}}#2\par}
\def\pmb#1{\setbox0=\hbox{#1}
\kern-0.025em\copy0\kern-\wd0
\kern.05em\copy0\kern-\wd0
\kern-.025em\raise.0433em\box0}
\def\ca#1{{\cal #1}}

\def\Cal{\ca}

\def\frac#1#2{{#1\over#2}}

\def\text#1{\rm#1}

\outer\def\stmnt  #1. #2\par{\medbreak
\noindent{\bf#1.\enspace}{\sl#2}\par
\ifdim\lastskip<\medskipamount \removelastskip\penalty55\medskip\fi}

\def\newline{\hfill\break}
\def\:{\,:\,}

\def\({\left(}
\def\){\right)}

\def\[{\left[}
\def\]{\right]}

\def\ci{\subset}

\def\fy{\infty}

\def\del{\partial}

\def\lam{\lambda}

\def\ph{\phi}

\cl{\bf On the  Mean Curvature Flow for $\sigma_k$-Convex
Hypersurfaces }

\ms
\cl{Hao Fang$^1$}
\cl{Department of Mathematics, Princeton University}
\cl{Princeton, NJ 08544, USA}
\ms
\cl{ChangYou Wang$^2$}
\cl{Department of Mathematics, University of Kentucky}
\cl{Lexington, KY 40506, USA}

\nind \S1 Introduction

For $n\ge 1$, let $M^n$ be a compact $n$-dimensional
manifold without boundary and $F_0: M^n\to R^{n+1}$
be a smooth immersion of $M^n$ in $R^{n+1}$ as a hypersurface.
Recall that $M_0=F_0(M^n)$ is said to be moved by its mean
curvature, if there is a family $F(\cdot,t)$ of
smooth immersions of $M^n$ into $R^{n+1}$ with the corresponding
hypersurfaces
$M_t=F(\cdot,t)(M^n)$ satisfying
$$\eqalignno{{\del F\over\del t}(p,t)&=-H(p,t)\nu(p,t),
\ \ (p,t)\in M^n\times R_+, &(1.1)\cr
                   F(p,0)&=F_0(p), \ \ \ \ \ \ \ \ \ \ \ \ \ \ \
\ \ \ \ p\in M^n &(1.2)\cr}$$
where $H(p,t)$ and $\nu(p,t)$ are the mean curvature and the outward unitary vector
of the hypersurface $M_t$ at $F(p,t)$, respectively. It was proved by
Huisken [H1] that there exists a $0<T=T(M^n)\le \fy$
such that (1.1)-(1.2) always admits a unique smooth solution for
$0<t<T$ and $\lim_{t\uparrow T}\max_{M_t}|H|=\fy$.
When $M_0$ is a convex hypersurface, Huisken [H1] proved that
$M_t$ contracts smoothly to a point as $t\uparrow T$.
Without the convexity assumption on $M_0$, Huisken proved
in [H2] (Theorem 3.5) that if the singularity is of the Type I, i.e.
$$\max_{M_t}|A|^2\le {C\over 2(T-t)}, \eqno(1.3)$$
then suitable scalings of $M_t$ near the singularity converges smoothly
to an immersed nonempty limiting hypersurface ${\tilde M}$, which
satisfies the equation
$$H(x)=<x,\nu(x)>, \eqno(1.4)$$
where $x$ is the position vector, $H$ is the mean curvature and
$\nu(x)$ is the outward unit normal vector field. It was also
proved by [H2] (Theorem 4.1) that compact manifolds with nonnegative
mean curvature satisfying (1.4) are spheres of radius $\sqrt n$.
For singularity of Type II, i.e.
$$\lim_{t\uparrow T}(T-t)\max_{M_t}|A|^2=\fy. \eqno(1.5)$$
Huisken and Sinestrari proved in [HS1] (Theorem 3.1) and
[HS2] (theorem 1.1) that if
$M_0$ is mean convex (i.e. $H\ge 0$) then $\sigma_k$-curvature ($2\le k\le
n$) of $M_t$ satisfies, for any $0<t<T$,
$$\sigma_k(M_t(p))\ge -\epsilon H^k(p,t)-C_{\epsilon,n,k}
\eqno(1.6)$$
(see below for the definition of $\sigma_k$). Based on this
key estimate, they proved in [HS2] (Theorem 4.1) that
the so-called essential scaling
of $M_t$ near time $T$ converges
to a smooth mean curvature flow $\{\tilde {M}_t\}_{t\in R}$,
with $\tilde{M}_t$ convex hypersurfaces. Moreover, either $\tilde{M}_t$
is strictly convex translating soliton or
(up to rigid motion) $\tilde{M}_t
=R^{n-k}\times \Sigma_t^k$ where $\Sigma_t^k$ is a lower
dimensional strictly convex soliton in $R^{k+1}$.
For $n=2$, it was shown by [HS1] (Corollary 4.7) that
$\{\Sigma_t^1\}_{t\in R}$ is the
``grim reaper'' curve given by $x=-\ln\cos y+t$.

In contrast with the convex case, it is well known
that the mean curvature flow (1.1)-(1.2) can develop
singularities before it may shrink to a point. It is a
major problem for people to study the nature of its
singularity and the asymptotic behavior near the
singularity. In this note, we make some effort to
try to understand the structure of the singularity
set at the first singular time for the initial hypersurface
$M_0$ belonging to the class consisting of $\sigma_k$ convex
hypersurfaces for some $1\le k\le n$.
More precisely, we want to understand
the limiting set $M_T$, which is the support of the
Radon measure $\mu_T$ obtained as the limit of Radon
measures $\mu_t$, where $\mu_t$ is the area measure of $M_t$ described
as below, as $t\uparrow T$. To better describe our result,
we first recall that, in addition to the above classical motion by
mean curvatures, Brakke [B] has introduced the motion by its mean
curvature for a family of Radon measures $\{\nu_t\}_{t\in R}$ in $R^{n+1}$
(e.g. $n$-dimensional rectifiable varifolds), which satisfies (1.1) in the
weak form
$$\int \phi\,d\nu_{t}\le \int\phi\,d\nu_{s}
+\int_s^t\int(-\phi |{\Cal H}|^2+{\Cal H}\cdot S^{\perp}\cdot
D\phi)\,d\nu_t\,dt
,\eqno(1.7)$$
for all nonnegative $\ph\in C_0^1(R^{n+1})$ and $0\le s\le t$.
The reader can refer to Ilmanen [I1] for the interpretation
of (1.7) and related results.

Note that if $\{M_t\}_{0\le t<T}$ solve (1.1)-(1.2) and
we denote $\mu_t$ as the area measure of $M_t$ for
$0\le t<T$, then $\mu_t$ are integral varifolds of multiplicity 1 and
satisfy
$$\mu_t(\phi)-\mu_s(\phi)=\int_s^t\int(-H^2\phi+D\phi\cdot {\Cal
H})\,d\mu_t\,dt, \eqno(1.8)$$
for any nonnegative $\phi\in C_0^1(R^{n+1})$ and $0\le s\le t<T$, where
$\phi_t(\phi)=\int \phi\,d\mu_t$ and ${\Cal H}=-H\nu$. In particular,
$\{\mu_t\}_{0\le t<T}$
is a motion by mean curvature by Brakke in the sense of (1.7). Moreover,
(1.8) implies that there exists a nonnegative measure $\mu_T$ in $R^{n+1}$
such that $\mu_t\rightarrow\mu_T$ as convergence of Radon measures in
$R^{n+1}$ as $t\uparrow T$. Now we extend $\{\mu_t\}_{0\le
t\le T}$ to $t>T$ such that $\{\mu_t\}_{0\le t\le \fy}$
are a family of Radon measures moved by its mean curvatures in the sense
of (1.7), whose existence is established by Brakke [B]. Our result is
concerned with the properties of
$M_T=\hbox{spt}(\mu_T)$.

One of the most important facts to the family $\{\mu_t\}_{0\le t\le \fy}$
is the following ${\it monotonicity \ formula}$, which was first
discovered by Huisken [H2] for (1.1) and later obtained by Ilmanen [I2]
for Brakke flows (1.7) and has played a  key role in
the analysis of the singularity for (1.1) and (1.7). The formula says
the follows. Let $\rho_{(y,s)}$ denote the $n$-dimensional backward
heat kernel centered at $(y,s)\in R^{n+1}\times R$ defined by
$$\rho_{(y,s)}(x,t)=(4\pi(s-t))^{-{n\over 2}}\exp(-{|x-y|^2\over 4(s-t)}),
\ \ x\in R^{n+1}, t<s.$$
Let $\mu=\{(\mu_t,t): 0\le t\le \fy\}$ and define, for $0<r<s$,
$$\Theta(\mu,(y,s),r)=\int \rho_{(y,s)}(x,s-r)\,d\mu_{s-r}(x).$$
Then one has
$$\eqalignno{&\int_{s-r_2}^{s-r_1}\int\rho_{(y,s)}(x,t)|{\Cal
H}(x,t)+{1\over 2(s-t)}S^{\perp}(x,t)\cdot (x-y)|^2\,d\mu_t(x)\,dt\cr
&\le\Theta(\mu,(y,s),r_2)-\Theta(\mu,(y,s),r_1). &(1.9)\cr}$$
for any $(y,s)\in R^{n+1}\times R_+$ and $0<r_1<r_2<s$, where
$S(x,t)=T_x\mu_t$.
A direct consequence of (1.9) is that
the density function
$$\Theta(\mu,(y,s))\equiv\lim_{r\downarrow 0}\Theta(\mu,(y,s),r)$$ exists
for any $(y,s)\in R^{n+1}\times R_+$ and is upper semicontinuous.
Using the upper semicontinuity, it is not difficult to prove
that $M_T$ is actually the Hausdorff limit of $M_t$ as $t\uparrow T$
(see, Lemma 2.1 below). Our first result is
\ss
\nind{\bf Theorem A}. {\sl The
$n$-dimensional Hausdorff measure of $M_T$ is finite,
i.e. $H^n(M_T)<\fy$.}

By exploring the upper semicontinuity of $\Theta(\mu,\cdot)$ and
extending the idea of the Federer's dimension reduction argument
to the parabolic setting, White [W1] (Theorem 9) has recently obtained the
stratification theorem for the support ${\Cal
M}=\{(\hbox{spt}(\nu_t),t):t\ge 0\}$
of Brakke flows for $k$-dimensional integral varifolds $\{\nu_t\}_{t\ge
0}$ in $R^{n+1}$ ($1\le k\le n$), which roughly says that the
points of ${\Cal M}$, for which each tangent flow having its
spine dimension (see [W1] for its definition) at most $l$, is of parabolic
Hausdorff dimension at most $l$ for all $0\le l\le k+2$.
Inspired by this stratification theorem by White [W1], we shall
consider the stratification of the extension set $M_T$ of the
$\{\mu_t\}_{t\ge 0}$ described as above. First, note that the monotonicity
of $\Theta(\mu,(y,s),\cdot)$ actually implies that $\Theta(\mu,(y,s))$ is
upper semicontinuous with both of its arguments
(see, e.g. [W2] theorem 2). Moreover, the
uniform upper bound of $\Theta(\mu,(y,s),\cdot)$ in terms of
$M_0$ and $(y,s)\in \hbox{spt}(\mu)$ implies that for $x\in M_T$
if we consider the parabolic blow-up, $P_{(x,T),\lam}(\mu)$, defined
as
$$P_{(x,T),\lam}(\mu)(\phi)=\lam^{-n}\int
\phi((x,T)+(\lam x,\lam^2 t))\,d\mu_t\,dt, \forall \phi\in
C_0^1(R^{n+1}\times R).$$
Then for any $\lam\rightarrow 0$ we can extract a subsequence
$\lam_i\rightarrow 0$ and a limiting Brakke flow
$\tilde{\mu}\equiv \{\tilde{\mu}_t\}_{t\in R}$, which is called
a tangent flow at $(x,T)$, such that
$P_{(x,T),\lam_i}(\mu)\rightarrow \tilde{\mu}$ as convergence
of Radon measures in $R^{n+1}\times R$. Moreover, as shown by
[I1] and [W1] that $\tilde{\mu}$ is backwardly self-similar, i.e.,
$P_{(0,0),\lam}(\tilde{\mu}|_{t\le 0})=\tilde{\mu}|_{t\le 0}$, and
$$\Theta(\tilde{\mu},(0,0))=\Theta(\mu,(x,T))\ge \Theta(\tilde{\mu},z),
\forall z=(y,s)\in R^{n+1}\times R.$$
It was also proved by [W1] that the set
$V(\tilde{\mu})\equiv \{x\in
R^{n+1}: \Theta(\tilde{\mu},(x,0))=\Theta(\tilde{\mu},(0,0))\}$
is a vector subspace of $R^{n+1}$. Denote
$\hbox{dim}({\tilde\mu})=\hbox{dim}(V(\tilde\mu))$. Then, White's
stratification
theorem yields
\ss
\nind{\bf Proposition B} ([W1]). {\it
Assume that the Brakke flow $\{\mu_t\}_{t\in R_+}$ is given
as above (i.e. $\{\mu_t\}_{0\le t<T}$ coincides
with the smooth mean curvature flow $\{M_t\}_{0\le t<T}$ of
(1.1)-(1.2)). Then
$$M_T=\Sigma_0\cup(\Sigma_1\setminus\Sigma_0)
\cup\cdots(\Sigma_n\setminus\Sigma_{n-1}),\eqno(1.10)$$
where $$\Sigma_i\equiv\{x\in M_T: \hbox{dim}(\tilde{\mu})\le i,
\hbox{ for any tangent flow }\tilde\mu \hbox{ at }(x,T)\}$$
for $0\le i\le n$.
Moreover, dim$_H(\Sigma_i)\le i$, for $0\le i\le n$.}

 Note that the original proof of White [W1], which is in the
nature of parabolic type, can be modified to prove this
proposition, we would like to give a slightly different proof
of it, which is in the Euclidean nature, in \S3.

Now we turn our attention to the $\sigma_k$ convex case of
mean curvature flows. First, we recall the definition of
$\sigma_k$ curvatures for hypersurfaces in $R^{n+1}$ (see
also [HS2]).
\ss
\nind{\bf Definition 1.1}. For a closed hypersurface $M\ci R^{n+1}$,
let $-\fy<\lam_1(x)\le \lam_2(x)\le\cdots\le \lam_n(x)<\fy$
be the principal curvatures of $M$ at $x\in M$. For $1\le k\le n$,
the $\sigma_k$ curvatures of $M$ at $x$ is defined by
$$\sigma_k^M(x)=\sum_{1\le i_1<\cdots<i_k\le n}\lam_{i_1}(x)\cdots
\lam_{i_k}(x).$$
Note that the $\sigma_1^M$ curvature is nothing but the mean curvature
of $M$, $\sigma_2^M$ is the scalar curvature of $M$, and $\sigma_n^M$
is the Guassian curvature of $M$.
\ss
\nind{\bf Definition 1.2}. For $1\le k\le n$. A closed hypersurface
$M\ci R^{n+1}$ is called $\sigma_k$ convex (respectively, strictly
$\sigma_k$ convex) if
$\hbox{min}_{x\in M}\sigma_k^M(x)\ge 0$ (respectively,
$\hbox{min}_{x\in M}\sigma_k^M(x)>0$). In particular,
the mean convexity of $M$ is equivalent to the $\sigma_1$
convexity of $M$.

For mean curvature flows of mean curvature sets $\{F_t(K)\}_{t\ge 0}$
(e.g. $K_0=F_0(\del K)$ is a mean convex smooth hypersurface), a
striking and difficult theorem by White [W3] (theorem 1) claimed that
the singular set of ${\Cal K}=\{(x,t): x\in
F_t(\del K),t\ge 0\}$, $\hbox{sing}
({\Cal K})$, has parabolic Hausdorff dimension at most $n-1$. Here
the singular set is defined to the completment of these regular points
near where ${\Cal K}$ is a smooth manifold and has its tangent plane
non-horizontal. As a direct consequence of this regularity theorem of
White, one knows that the top dimensional subset
$\Sigma_n\setminus\Sigma_{n-1}\ci M_T$ defined as in
the proposition B is regular set of $M_T$, namely
near each point in $\Sigma_n\setminus\Sigma_{n-1}$ $M_T$ is
a smooth $n$-dimensional manifold.

 By exploring the estimates (1.6) on $\sigma_k$ curvatures for the
mean curvature flow of mean convex hypersurfaces (1.1)-(1.2) obtained
by [HS2] and the partial regularity theorem for mean convex flows of [W3]
(i.e. the singular set has parabolic Hausdorff dimension at most $n-1$),
we obtain the following result.
\ss
\nind{\bf Theorem C}. {\sl Assume that the Brakke flow $\{u_t\}_{t\in
R_+}$ is given as same as that in Proposition B. For any $2\le k\le n-1$,
if the initial closed hypersurface $M_0$ is $\sigma_k$ convex.
Then
$$\Sigma_{n-1}\setminus\Sigma_{n-2}=\cdots=\Sigma_{n-k+1}\setminus
\Sigma_{n-k}=\emptyset. \eqno(1.11)$$}

 This note is written as follows. In \S2, we prove Theorem A;
In \S3, we give a proof of Proposition B; In \S4, we prove Theorem C.

\bs
\nind{\S2} Proof of Theorem A

This section is devoted to the proof of theorem A. First, we show
\ss
\nind{\bf Lemma 2.1}. {\sl $M_t$ converges to $M_T$ in the Hausdorff
distance sense, as $t\uparrow T$. }
\ss
\nind{\bf Proof}. For any $t\uparrow T$, we can extract a subsequence
$t_i\uparrow T$ and a closed subset $A\ci R^{n+1}$ such that
$M_{t_i}$ converges to $A$ in the Hausdorff distance.
Now we want to show $A=M_T$.
Suppose $x_0\not\in A$. Then there exists $r_0>0$ such that
$B_{r_0}(x_0)\cap M_{t_i}=\emptyset$ for $i$ sufficiently large.
In particular, $H^n(M_{t_i}\cap B_{r_0}(x_0))=0$.
Hence, $\mu_T(B_{r_0}(x_0))=0$ and
$x_0\not\in M_T$. This gives that $M_T\ci A$.
To prove $A\ci M_T$, we argue
by contradiction. Suppose that there exists $x_0\in A\setminus M_T$.
Then there exists
$r_0>0$ such that $\mu_T(B_{r_0}(x_0))=0$. On the other hand, there exist
$x_i\in M_{t_i}$ such that $x_i\rightarrow x_0$. Therefore,
by the upper semicontinuity
of $\Theta(\mu, \cdot)$ and the fact that
$\Theta(\mu,(x_i,t_i))=1$ (since $x_i\in M_{t_i}$ and $M_{t_i}$
is smooth), we have
$$\Theta(\mu, (x_0,T))\ge \limsup_{i\rightarrow
\fy}\Theta(\mu, (x_i,t_i))=1.$$
This implies that $x_0\in M_T$ (for otherwise $\Theta(\mu, (x_0,T))=0$).
We get the desired contradiction.  \qed
\ss
\nind{\bf Proof of Theorem A}. Since
for any $0\le t<T$ and $y\in M_t$, one
has $(\Theta,(y,t))=1$. The upper semicontinuity implies that
$\Theta(\mu,(x_0,T))\ge 1$ for any $x_0\in M_T$. For any $\epsilon>0$,
the estimate of Cheng [C] implies
that there exists sufficiently large $K_\epsilon>0$
such that
$$\eqalignno{G_{(x_0,T)}(x,T-r^2)&\le r^{-n}, \ \ \ \ \ \ \ \ \ \ \
\  \ \ \ \ \ \ \ \ \ \ \ \forall x\in R^{n+1},\cr
 &\le \epsilon G_{(x_0, T+r^2)}(x,T-r^2),
\forall |x-x_0|\ge K_\epsilon r.\cr}$$
Therefore, we have, for any $x_0\in M_T$,
$$\eqalignno{1&\le \Theta(\mu, (x_0,T))\cr
&\le Cr^{-n}H^n(M_{T-r^2}\cap
B_{K_\epsilon r}(x_0))\cr
&+\epsilon\int_{M_{T-r^2}}G_{(x_0,T+r^2)}(x,T-r^2). &(2.1)\cr}$$
The monotonicity for $\Theta(\mu,(x_0,T+r^2),\cdot)$ implies
that
$$\eqalignno{&\int_{M_{T-r^2}}G_{(x_0,T+r^2)}(x,T-r^2)\cr
&=\int_{M_{T+r^2-2r^2}}G_{(x_0,T+r^2)}(x,T+r^2-2r^2)\cr
&\le\int_{M_{T+r^2-R^2}}G_{(x_0.T+r^2)}(x,T+r^2-R^2)\cr
&\le CR^{-n}H^n(M_{T+r^2-R^2})\le CR^{-n}H^n(M_0)\cr}$$
for any $\sqrt{2}r\le R\le \sqrt{T+r^2}$.
In particular, for sufficiently small $r$, by choosing $R^2={T\over 2}$,
we have
$$\int_{M_{T-r^2}}G_{(x_0,T+r^2)}(x,T-r^2)\le
CT^{-{n\over 2}}H^n(M_0).$$
Hence, by choosing $\epsilon=\epsilon(M_0,T)>0$ sufficiently small, we
have
$$ r^n\le CH^n(M_{T-r^2}\cap B_{K_\epsilon r}(x_0)). \eqno(2.2)$$
Observe that the family
${\Cal F}=\{B_{K_\epsilon r}(x): x\in M_T\}$ covers $M_T$ so that
the Vitali's covering Lemma implies that there exists a disjoint
subfamily $\{B_{K_\epsilon r}(x_i): x_i\in M_T\}_{i=1}^\fy$ such that
$$M_T\ci\cup_{i} B_{5K_\epsilon r}(x_i),$$
so that
$$\eqalignno{H_{5K_\epsilon r}^n(M_T)
&\le (5K_\epsilon)^n\sum_{i} r^n\cr
&\le (5K_\epsilon)^n\sum_{i}H^n(M_{T-r^2}\cap B_{K_\epsilon r}(x_i))\cr
&= (5K_\epsilon)^n H^n(M_{T-r^2}\cap (\cup_{i} B_{K_\epsilon r}(x_i)))\cr
&\le (5K_\epsilon)^n H^n(M_{T-r^2})
\le C(\epsilon, M_0)<\fy.\cr}$$
This finishes the proof for Theorem A.  \qed
\bs
\nind{\S3} Proof of Proposition B

 In this section, we give a slightly different and also
easier proof of Proposition B, which is essentially due to White [W1].

 For any $0\le i\le n$, it follows from the definition of $\Sigma_i$
that for any $x_0\in\Sigma_i$ and each tangent flow $\tilde\mu$ of
$\mu$ at $(x_0,T)$, which is a backwardly self-similar Brakke flow,
$V(\tilde\mu)=\{x\in R^{n+1}:\Theta(\tilde\mu,(x,0))
=\Theta(\tilde\mu,(0,0))\}$ is a vector space of dimension at most
$i$. Moreover, it follows from the argument in [W1] (see also
theorem 3 [W2]) that if we let
$W(\tilde\mu)\equiv\{(x,t)\in R^{n+1}\times R_{-}:\Theta(\tilde\mu,(x,t))
=\Theta(\tilde\mu,(0,0))\}$ then either
$W(\tilde\mu)=V(\tilde\mu)\times R_{-}$ and
there exists a minimal hypercone $C^{n-d}\in R^{n+1-d}$
such that
$$\tilde\mu|_{t\in R_{-}}=(R^d\times C^{n-d})\times R_{-},\eqno(3.1)$$
or $W(\tilde\mu)=V(\tilde\mu)$ and there exists a backwardly
self-similar $(n-d)$ dimensional Brakke flow $\nu=\{\nu_t\}_{t\in R}$
in $R^{n+1-d}$ such that
$${\tilde\mu}|_{t\in R_{-}}=V({\tilde\mu})\times \nu|_{t\in R_{-}},
\eqno(3.2)$$
here $d=\hbox{dim}(V(\tilde\mu))$. Define
$$\eta_{x,\rho}(y)=\rho^{-1}(y-x), \forall y\in R^{n+1}$$
Now we claim
\ss
\nind{\bf Claim 3.1}. {\sl For any $x_0\in \Sigma_i$ and each $\delta>0$
there exists an $\epsilon>0$ (depending on ${\mu}$, $x_0$, $\delta$)
such that for each $\rho\in (0,\epsilon]$
$$\eqalignno{&\eta_{x_0,\rho}\{x\in
B_\rho(x_0):\Theta(\mu,(x,T))\ge\Theta(\mu,(x_0,T))-\epsilon\}\cr
&\ci\hbox{ the }\delta-\hbox{neighbourhood of }L_{x_0,\rho} &(3.3)\cr}$$
for some $i$-dimensional subspace $L_{x_0,\rho}$ of $R^{n+1}$.}
\ss
\nind{\bf Proof}. If this is false, there exist $\delta>0$ and $x_0\in
\Sigma_i$ and $\rho_k\downarrow 0$ and $\epsilon_k\downarrow 0$ such that
$$\eqalignno{&\{x\in B_1(0):\Theta(P_{(x_0,T),\rho_k}(\mu),(x,0))
\ge\Theta(\mu,(x_0,T))-\epsilon_k\}\cr
&\not\ci \delta-\hbox{neighbourhood
of }L, &(3.4)\cr}$$
for any $i$-dimensional subspace $L$ of $R^{n+1}$. But
$P_{(x_0,T),\rho_k}(\mu)\rightarrow \tilde\mu$, a tangent flow
of $\mu$ at $(x_0,T)$, and $\Theta(\tilde\mu,(0,0))=\Theta(\mu,(x_0,T))$.
Since $x_0\in\Sigma_i$, we have $\hbox{dim}(V(\tilde\mu))\le i$,
there is a $i$-dimensional subspace $L_0\ci R^{n+1}$ such that
$V(\tilde\mu)\ci L_0$. Moreover, the uppersemicontinuity of
$\Theta(\tilde\mu,\cdot)$ implies that there is a $\alpha>0$
such that
$$\Theta(\tilde\mu,(x,0))<\Theta(\tilde\mu,(0,0))-\alpha,
\forall x\in B_1(0) \hbox{with dist}(x,L_0)\ge\delta. \eqno(3.5)$$
Then the upper semicontinuity of $\Theta(\mu,\cdot)$ for convergence
of both of its variables
implies that
we must have, for $k'$ sufficiently large,
$$\Theta(P_{(x_0,T),\rho_{k'}}(\mu), (x,0))<\Theta(\tilde\mu, (0,0))
-\alpha, \forall x\in B_1(0) \hbox{ with dist}(x,L_0)\ge\delta.
\eqno(3.6)$$
This clearly contradicts with (3.4). The claim is proven.  \qed
\ss
\nind{\bf Completion of Proof of Proposition B}. We decompose
$\Sigma_i=\cup_{l=1}^\fy\Sigma_i^l$, where $\Sigma_i^l$ denotes
the points $x\in\Sigma_i$ such that the claim 3.1 holds for
$\epsilon=l^{-1}$. Now we decompose
$\Sigma_i^l=\cup_{q=1}^\fy\Sigma_i^{l,q}$, where
$\Sigma_i^{l,q}=\{x\in\Sigma_i^l:
{q-1\over i}\le\Theta(\mu,x)\le {q\over i}\}$. Hence, claim 3.1
implies that for $A=\Sigma_i^{l,q}$,
$$\eta_{x,\rho}(A\cap B_\rho(x))\ci
\delta-\hbox{neighbourhood of } L_{x,\rho}, \forall x\in A,
\rho<l^{-1}. \eqno(3.7)$$
for some $i$-dimensional subspace $L_{x,\rho}\ci R^{n+1}$.
The proof is completed if we apply the following Lemma, whose
proof can be found in the Lecture 2.4 of Simon [S].
\ss
\nind{\bf Lemma 3.2}. {\sl There is a $\beta:R_+\to R_+$ with
$\lim_{t\downarrow 0}\beta(t)=0$ such that
if $\delta>0$ and if $A\in R^{n+1}$ satisfying the
property (3.7) above, then $H^{i+\beta(\delta)}(A)=0$.}

\qed

\bs
\nind{\S4} Proof of Theorem C

In this section, we outline the proof of the theorem C. But,
first, we gather together needed key estimates by [HS2]
on the $\sigma_k$ curvature under the mean curvature flow
(1.1)-(1.2) with the initial $M_0$ being mean convex.
\ss
\nind{\bf Lemma 4.1}. {\sl (a). For any closed hypersurface $M$
in $R^{n+1}$. For any $1\le k\in n$, if $M$ is $\sigma_k$ convex (or
$\sigma_k$ strictly convex, respectively)
then $M$ is also $\sigma_l$ convex (or $\sigma_l$ strictly
convex, respectively) for all $1\le l\le k$. In particular, $M$ is mean
convex.}

{\sl \nind (b). For any $1\le k\le n$. Assume that $M_0$ is a $\sigma_k$
convex closed hypersurface and $\{M_t\}_{0\le t<T}$ is the
mean curvature flow (1.1)-(1.2). Then, for any $0<t<T$, $M_t$
is $\sigma_k$ strictly convex. In particular, $M_t$ is $\sigma_l$
strictly convex for $0<t<T$ and all $1\le l\le k$.}
\ss
\nind{\bf Proof}. The reader can find the details of the proof
in Proposition 3.3 (i) (ii) in [HS2].   \qed

It follows from Lemma Lemma 4.1 that for $1\le k\le n$ if $M_0$ is a
$\sigma_k$ convex  closed hypersurface then for any $0<t<T$ there exists a
$\epsilon=\epsilon_t>0$ such that
$$\sigma_l^{M_t}(x)\ge \epsilon H^l(x,t), \forall x\in M_t, 1\le l\le k.
\eqno(4.1)$$
Proposition 3.4 of [HS2] then asserts that the inequality (4.1)
is kept under the mean curvature flow (1.1)-(1.2). More
precisely, we have
\ss
\nind{\bf Lemma 4.2}. {\sl For $1\le k\le n$.
Let $\{M_t\}_{0\le t<T}$ be the mean curvature flow (1.1)-(1.2),
with $M_0$ being a $\sigma_k$ convex closed hypersurface.
Then (3.1) holds with $\epsilon=\epsilon_t>0$ for all $s\in [t,T)$
and all $1\le l\le k$.}
\ss
\nind{\bf Proof of Theorem C}. Suppose that the conclusion fails.
Then there exists a $2\le i\le k$ such that
$\Sigma_{n-i+1}\setminus\Sigma_{n-i}\not=\emptyset$.
Pick a point $x_0\in\Sigma_{n-i+1}\setminus \Sigma_{n-i}$.
Then it follows from the definition that any tangent flow of
$\mu$ at $(x_0,T)$ has no more than $n-i+1$ dimension
of spatial translating invariant directions, and there
exists at least one $\lam_m\downarrow 0$ and a vector subspace $V\ci
R^{n+1}$, with $\hbox{dim}(V)=n-i+1$, such that  either (3.1) or (3.2)
holds, namely either (a): there exists a minimal hypercone
$C^{i-1}\ci R^{i}$ such that
$$P_{(x_0,T),\lam_m}(\mu)|_{t\in R_{-}}\rightarrow
(V\times C^{i-1})\times R_{-}, \eqno(4.1)$$
or (b): there exists a backwardly self-similar $(i-1)$ dimensional
Brakke flow $\nu=\{\nu_t\}_{t\in R}$ in $R^i$ such that
$$P_{(x_0,T),\lam_m}(\mu)|_{t\in R_{-}}\rightarrow
V\times \nu|_{t\in R_{-}}. \eqno(4.2)$$
As a special case of the partial regularity theorem of White [W1],
we know that the singular set of both $C^{i-1}$ (in the case (a)) and
${\Cal M}_{-1}
=\hbox{spt}(V\times \nu_{-1})$ (in the case (b)) has Hausdorff dimension
at most
$i-4$, here the singular set of a subset $A\in R^i$ is the completment
of regular points in $A$ and the regular points are points near where
$A$ is a smooth manifold. For, otherwise, the singular set
of the corresponding tangent flow has parabolic Hausdorff
dimension larger than $(i-4)+(n-i+1)+2=n-1$, which
contradicts with White's theorem. Since case (a) can be handled
similarly to case (b), we want to discuss case (b) only.
Let $N\ci {\Cal M_{-1}}$ denotes the singular set. Then, for
any $x\in {\Cal M_{-1}}\setminus N$, $V\times\nu$ is smooth
in a spacetime neighbourhood $U$ of $(x,-1)$. Moreover, by
the Brakke (unit density) regularity theorem ([B]), we can assume
that $\mu_m\equiv P_{(x_0,T),\lam_m}(\mu)$ converge smoothly to $V\times
\nu$
in the neighbourhood $U$. We now claim
\ss
\nind{\bf Claim 4.3}. {\sl The second fundamental form $A$ of
${\Cal M}_{-1}$ vanishes
everywhere on ${\Cal M}_{-1}\setminus N$}.
\ss
\nind{\bf Proof of Claim 4.3}. Since $M_0$ is assumed to be $\sigma_k$
convex, it follows from part (b) of Lemma 4.1 that without loss
of generality we can assume that $M_0$ is in fact $\sigma_k$
strictly convex so that Lemma 4.2 implies that there exists
a $\epsilon>0$ such that
$$\sigma_l^{M_t}(x)\ge\epsilon H^l(x,t)>0, \forall x\in M_t,0\le t<T,
 1\le
l\le k.
\eqno(4.3)$$
Note that (4.3), combined with scalings, implies that,
for any $1\le l\le k$, $t\in [-\lam_m^{-2}T_0, 0)$,
and $x\in \hbox{spt}(\mu_m(t))=\hbox{spt}(\mu_m\cap\{t\})$,
$$\eqalignno{\sigma_l^{\mu_m(t)}(x)&=\lam_m^l\sigma_l^{M_{T_0+\lam_m^2
t}}(x_0+\lam_m x)\cr
&\ge \epsilon \lam_m^l H^l(x_0+\lam_m x, T_0+\lam_m^2 t)\cr
&=\epsilon (\sigma_1^{\mu_m(t)}(x))^l, &(4.4)\cr}$$
$$H(x,t)>0, \forall x\in \hbox{spt}(\mu_m(t)), t\in [\lam_m^{-2}T_0,0).
\eqno(4.5)$$
This, combined with the smooth convergence fact at
${\Cal M}_{-1}\setminus N$ as above,
implies, for any $1\le l\le k$
$$\sigma_l^{{\Cal M}_{-1}}(x)\ge \epsilon H^l(x,-1)\ge 0,
 \forall x\in {\Cal M_{-1}}\setminus N. \eqno(4.6)$$
On the other hand, since ${\Cal M}_{-1}=V\times \hbox{spt}(\nu_{-1})$
and $\hbox{dim}(V)=n-i+1$, it follows from the definition 1.1 that
$$\sigma_i^{{\Cal M}_{-1}}(x)=\cdots=\sigma_k^{{\Cal M}_{-1}}(x)
=0, \forall x\in {\Cal M}_{-1}\setminus N. \eqno(4.7)$$
(4.6) and (4.7) imply that
$$H(x,-1)=0, \forall x\in {\Cal M}_{-1}\setminus N. \eqno(4.8)$$
Therefore, for any $x\in {\Cal M}_{-1}\setminus N$,
$$\eqalignno{\sigma_2^{{\Cal M}_{-1}}(x)
&={1\over 2}(H^2(x,-1)-|A|^2(x,-1))\cr
&=-{1\over 2}|A|^2(x,-1)\le 0.\cr}$$
This, combined with (4.6),  implies that
$$\sigma_2^{{\Cal M}_{-1}}(x)=|A|(x,-1)=0, \forall x\in {\Cal
M}_{-1}\setminus N. \eqno(4.9)$$
This finishes the proof of Claim 4.3.
It follows from the Claim 4.3 that ${\Cal M}_{-1}\setminus N$ is
flat. Note also that $V\times \nu_t|_{t\in R_-}$ is also
a backwardly self-similar Brakke flow. Therefore, we know
that ${\Cal M}_{-1}$ is also a minimal hypercone in $R^{n+1}$.
Hence ${\Cal M}_{-1}$ is a hyperplane in $R^{n+1}$. This
clearly contradicts with the fact that $2\le i\le k$ and $x_0\in
\Sigma_{n-i+1}\setminus\Sigma_{n-i}$. Therefore, the proof is
complete.  \qed

\ms
\nind {\bf ACKNOWLEDGEMENTS}.  Both authors learned this subject
through a topic course by Prof. G. Huisken  during the Spring
semester, 2000. We would like to thank him for 
his encouragements to us for writing this note.

\bs
\cl{\bf REFERENCES}
\ss
\nind{[B]} K. Brakke, {The motion of a surface by its mean curvature}.
Princeton U. Press, 1978.

\ss\nind{[C]} X. Cheng, {\it Estimate of the singular
set of the evolution problem for harmonic maps}. J. Diff. Geom.
34 (1991), no. 1, 169-174.

\ss\nind{[H1]} G. Huisken, {\it Flow by mean curvature of convex surface
into spheres}. J. Diff. Geom., 20 (1984), 237-266.

\ss\nind{[H2]} G. Huisken, {\it Asymptotic behaviour for singularities
of the mean curvature flow}. J. Diff. Geom. 31 (1990), 285-299.

\ss\nind{[HS1]} G. Huisken \& C. Sinestrari, {\it Mean curvature flow
singularities for mean convex surfaces}. Calc. Var. Partial Differential
Equations, 8 (1999), 1-14.

\ss\nind{[HS2]} G. Huisken \& C. Sinestrari, {\it Convexity estimates for
mean convex flow and singularities of mean convex surfaces.}

\ss\nind{[I1]} T. Ilmanen, {Elliptic regularization and partial regularity
for motion by mean curvature}. Memoris of AMS , 108 (1994).

\ss\nind{[I2]} T. Ilmanen, {\it Singularities of mean curvature flow
of surfaces}. Preprint.

\ss\nind{[S]} L. Simon, {\it Singularities of geometric variational
problems}. Nonlinear partial differential equations in differential
geometry (Park City, UT, 1992), 185-223. IAS/Park City
Math. Ser., 2, Amer. Math. Soc., Providence, RI, 1996.

\ss\nind{[W1]} B. White, {\it  Stratification of minimal surfaces, mean
curvature flows, and harmonic maps}. J. Reine angew. Math. 488 (1997),
1-35.

\ss\nind{[W2]} B. White, {\it  Partial regularity theorem for mean
curvature flows}. International Math. Research Notices 4 (1994), 185-192.

\ss\nind{[W3]} B. White, {\it The size of the singular set in mean
curvature flow of mean convex surfaces}. J. Amer. Math. Soc. 13 (2000),
no. 3, 665-695.

\end